\pdfoutput=1
\documentclass[10pt,psamsfonts,oneside]{amsart}
\usepackage[left=1.5in,right=1.5in,top=1in,bottom=1in]{geometry}
\usepackage{glossaries}
\begin{document}
\title{A Travel-Time Metric}
\author{Kenneth Halpern}
\address{Cambridge, MA}
\thanks{We thank A. Constandache for comments and suggestions.}
\begin{abstract}
Inspired by the question of whether one can create a visually informative map\footnote{We use the term ``map'' in the cartographic sense rather than the mathematical one.} of a geographic region, in which distances reflect travel-times rather than physical proximity, we examine whether it is possible to construct a meaningful travel time function between points in a region.  We argue that such a function should constitute a mathematical metric, and we offer a reasonable candidate using a max-min approach.
\end{abstract}

\maketitle

\parskip 10pt

\section{Introduction}

Often the physical proximity of two geographic locations is less interesting than the time it takes to travel between them.  The two can differ markedly, and often do.  It could be useful to define some sort of ``travel-time'' function between (ordered) pairs of locations in a region.   If chosen wisely, the latter could and should inform a wide variety of decisions.  These include personal choices such as where we live, work, vacation, and shop; business strategies for the planning of delivery routes, placement of retail or storage locations, and general expansion; and governmental policy or administration at various levels relating to the optimal placement of infrastructure, facilitation of traffic flow, urban planning, and accessibility.  It is not hard to imagine myriad other applications as well.  Except where distances are large and the cost of fuel a major consideration, it is fair to say that the travel time between points is of greater practical importance than their physical distance.  Though various non-Euclidean metrics (most notably ``taxicab'' geometry\cite{taxicab} and Chebyshev distance\cite{chebyshev}) have been explored, these merely offer alternate views of physical distance and are poor proxies for the travel-time.  

While individual travel-times between points find use in many present-day applications and are readily available from a variety of sources, an overall travel-time function would prove far more informative and provide a comprehensive basis for the types of decisions mentioned.   Ideally, a cartography based on it would be possible, allowing for a meaningful two-dimensional representation of the mutual travel-time proximities of points in some set\footnote{Strictly speaking, we cannot simply distort an arbitrary travel-time geometry into the equivalent of an ordinary map.  In our framework, this would amount to constructing an autohomeomorphism that transforms an arbitrary metric into some particular metric (typically quadratic), and that is impossible in general.  However, it is feasible under certain restrictive conditions.  We also may be able to devise other types of visualizations that are useful.}. That is beyond our present scope, however.  Here, we focus on the development of a travel-time metric on a set of points.  

In order for a travel-time function to be useful for general purposes, it must be independent of the particular route or mode of transportation chosen.  As we will see, there are strong reasons to desire it to be independent of departure time as well -- at least over short periods.  While predictive models may have value for navigation or other applications, we confine ourselves to the consideration of ex-post-facto travel-time functions.  These are based on realized travel-times over some sampling period\footnote{This period may roll or grow as needed, but at any given time of analysis there should be a fixed travel-time function.}.  For these reasons, the travel-time function will involve some sort of aggregate -- a term we use loosely to designate a function derived from a set of sample data in some fashion.  The aggregation will be over routes and departure times within the sample period.  

Notwithstanding these considerations, the variation of such a travel-time function by hour of the day, by season, or due to temporary conditions also may prove enlightening.  Growing problems could be detected and mitigated, investment opportunities perceived, and regional changes observed.  Only recently have we reached the point at which the relevant technology is sufficiently advanced and its adoption sufficiently widespread to provide measurement data for our approach.  Part of our purpose is to encourage the construction, dissemination, and use of travel-time functions.  

\section{Routes}

We consider a set of points $M$ between which travel is possible\footnote{The ordinary geographic topology and geometry of $M$ do not concern us here.}, and a set of departure times within some period $P$.

A ``route'' $r(a,b,t)$ between points $a,b\in M$ is defined to be a physical path between those points, along with the modes of transportation by which that journey is effected beginning at time $t$.  Note that in our notation, $a$, $b$, and $t$ are part of the label for the route, rather than functional parameters designating progress along it. The physical path may cross itself or involve repeated portions, and there may be multiple simultaneous routes along the same physical path between two points\footnote{These could represent different modes of transportation, for instance.}.  

We denote by $R(a,b,t)$ the set of all routes under consideration between $a$ and $b$ with departure time $t$.  As will be seen, there is flexibility in the routes we include in $R$ as long as we adhere to certain consistency requirements.  

We define $T(r(a,b,t))$ to be the travel time associated with route $r(a,b,t)$.   Much like that of $R$, the precise definition of $T(r(a,b,t))$ admits some flexibility\footnote{A predictive model could make various assumptions about the times associated with the travel components of $r$, for example.}.   However, in almost any real use case, $T(r(a,b,t))$ would be based on empirical data and offer no such freedom.  

Two successive routes that lay end-to-end may be joined to form a possible third route.  We write $(r+r')$ for the route $r(a,b,t)$ followed by the route $r'(b,c,t+T(r(a,b,t)))$. 

There are a few consistency requirements that we must impose for things to make sense.  

\begin{itemize}
\item \textbf{Identity}: There exists a $0$-time route $e(a,a,t)$ between any point $a\in M$ and itself for any departure time $t\in P$.  $T(e(a,a,t))=0$.
\item \textbf{Positivity}: All routes other than the $e(a,a,t)$ have travel time $T(r(a,b,t))>0$.  It takes some time to travel between two different points.  
\item \textbf{Existence}: Given any two points $a,b\in M$ and any departure time $t\in P$, there exists at least one finite travel-time route $r(a,b,t)$.  That is, we always can travel between any two points in either direction.
\item \textbf{Composition}:  Let $a,b,c\in M$ and $t\in P$.  If $r(a,b,t)\in R(a,b,t)$ and $t+T(r(a,b,t))\in P$ and $r'(b,c,t+T(r(a,b,t)))\in R(b,c,t+T(r(a,b,t)))$, then $(r+r')(a,c,t)\in R(a,c,t)$.  If we travel one route followed by another (and both the end time and end point of the first route coincide with the start time and start point of the second), then we may travel them in succession.  The overall travel time is the sum of the two component travel times.  Note that we only require composition if the departure times of \textit{both} the first and second routes are in $P$.  The second need not finish within it, though.  
\end{itemize}

Though unnecessary for any of our proofs, it also makes sense to allow a set of routes which correspond to waiting in place.  That is, we define a large set of routes $w(a,a,t,T)$ for all $a\in M$, $t\in P$, and $T>0$.  These can be composed with other routes to allow waiting in between legs of a journey.  

Mathematically, we have a strongly connected directed graph $G$ whose vertices are the points in $M$ and whose edges are the routes.  This is true overall, but also when we restrict ourselves to routes with any given departure time $t\in P$.  Unfortunately, our time-dependent composition rule does not readily translate into a simple structure on either the overall graph or the time-snapshot graphs.  For example, $G$ is not a category.

As will be demonstrated, we may choose the set of routes $R$ in any manner that adheres to our consistency requirements.  The utility of the resulting metric will depend on the prudence of our choice, of course.  

\section{Travel-Time Function at a given Time}

There are many candidates for a travel-time function $T(a,b)$, each with various benefits and drawbacks. We focus on whether it is possible to find a sensible choice that also is a metric.  Recall that $T$ is a metric on set $M$ if it is a function $M\times M\rightarrow R$ that satisfies four requirements:  (i) $T(a,b)\geq 0$, (ii) $T(a,b)=0$ iff $a=b$, (iii) $T(a,b)=T(b,a)$, and (iv) $T(a,c)\leq T(a,b)+T(b,c)$.  Why should we demand these particular properties for $T$?  Where needed, we turn to our original cartographic inspiration.  

The requirements that $T(a,a)=0$ and $T(a,b)> 0$ if $a\ne b$ are self-evident.  It takes no time to get from a place to itself, and it always takes time to get between two different places.  Any meaningful measure of travel-times must adhere to these.  The triangle inequality also makes intuitive sense.  Simplistically, if $T(a,b)+T(b,c)$ is faster than $T(a,c)$ we should define $T(a,c)$ using the better route via $b$.  There are many good reasons this may be violated when we aggregate, and our imposition of it is a choice.  That $T$ should be symmetric is the least justifiable of our requirements.  On a route-by-route basis, it clearly is wrong; and there is no reason it should hold in general.  An aggregate $T$ may or may not be symmetric, depending on the definition.  Like the triangle-inequality, symmetry is necessary for cartography.  There also are more general reasons we may wish $T$ to be a metric.  In order to meaningfully speak of distances between points -- particularly in the presence of composable routes -- it is necessary to have one.  However, these all are mere motivations; the decision to pursue a metric remains our choice. 

We begin by looking at the route-specific travel-time $T(r(a,b,t))$.  Our composition requirement states that for successive end-to-end routes $r$ and $r'$:

\begin{equation}
T(r(a,b,t))+T(r'(b,c,t+T(r(a,b,t))))= T((r+r')(a,c,t))
\label{eq:tti}
\end{equation}

If we hope to obtain a metric, we must aggregate these over routes and departure times in some fashion.  Minimizing over routes with a given departure time is a good prospective first step.  Other aggregates have a fatal flaw.  It is easy to construct arbitrarily long routes, and the set of routes with a given departure time is infinite\footnote{This follows from the existence and composition requirements.}.  Any statistical aggregate would involve a summation over routes, receive unbounded contributions, and prove divergent unless artificially regulated in some way.   A maximization over routes similarly would prove unbounded.  There may be other ways of aggregating, but minimization is the simplest and most obvious choice.  It also makes intuitive sense.  At any given departure time, one may choose from amongst a certain set of routes.  Presumably, the quickest would be selected.  While the absence of perfect information may lead to suboptimal choices, and perhaps justify a statistical approach, we already have seen that this is not easily accomplished without introducing additional arbitrary assumptions and a consequent dependence upon them.  

In this vein, it is convenient to define a best travel-time function $T(a,b,t)$ as\footnote{Unlike our notation for $r$ and $R$, the $a$, $b$, and $t$ here designate actual variables of which $T$ is a function.}

\begin{equation}
T(a,b,t)= \min_{r\in R(a,b,t)} T(r(a,b,t))
\label{eq:timedepdef}
\end{equation}

This neither is symmetric nor obeys the triangle inequality.  However, it is not hard to see that it obeys a time-dependent sibling of the latter.  If $t+T(a,b,t)\in P$ then

\begin{equation}
T(a,c,t)\leq T(a,b,t)+ T(b,c,t+T(a,b,t))
\label{eq:newtriang}
\end{equation}

Each term on the right is a minimization over routes starting at the relevant time and point.  We may select any route $r(a,b,t)$ which possesses the minimum travel-time $T(a,b,t)$ (there must be at least one) for the first leg and compose it with any route $r'(b,c,t+T(a,b,t))$ which possesses the minimum travel-time $T(b,c,t+T(a,b,t))$ for the second.  The result is a legitimate route from $a$ to $c$ departing at $t$.  Since $T(a,c,t)$ is a minimum over all such routes, it must include this one and cannot take longer than it.

\section{Travel-Time Function aggregated over Time}

Next we must remove the dependence on departure time.  Here, as well, statistical aggregates prove problematic.  Each such quantity corresponds to a weighted sum or integral over some function of the $T(a,b,t)$'s (possibly involving $T$'s with different departure times).  It is not hard to see that most such integrals can be made to violate the triangle inequality by a judicious choice of the $T$'s themselves.  For example, consider an obvious candidate

$$T_I(a,b)= \int_P w(t) T(a,b,t)dt$$

for some weight function $w(t)\ge 0$ with $\int w(t)=1$.

For simplicity, let's assume there is no boundary\footnote{We will need to deal with boundary issues even in the max-min approach we ultimately adopt -- so we should assume that we could deal with them sensibly here too if need be.} and that $P=[-\infty,\infty]$.  From equation~\ref{eq:newtriang} we obtain

$$T(a,c)\leq T(a,b) + \int T(b,c,t+T(a,b,t)) w(t) dt$$

In order to guarantee that the triangle inequality holds we require $\int T(b,c,t) w(t)\ge \int T(b,c,t+T(a,b,t)) w(t)$.  This must be true for all possible functions $T(*,*,t)$ which may arise from our consistency requirements.  The only consequent constraints on $T(*,*,t)$ are that it is $\ge 0$ and finite for all pairs of points and all departure times.  

However, it is not hard to engineer a $T$ for which that inequality fails.  If $b=c$ or $a=b$ then the Triangle Inequality holds trivially, so we need only consider distinct points.  Not only is it easy to violate the requirement that $\int T(b,c) w(t) \ge \int T(b,c,t+T(a,b,t)) w(t)$, but we can do so at each departure time distinctly.  The problem may be rephrased as follows: Given a function $g(t)>0$, we must find a function $f(t)>0$ for which $f(t+g(t))> f(t)$ for all $t$.  However, any monotonically increasing $f$ will do the trick.  Analogous problems plague any reasonable function of $T$ that could serve as an integrand as well.  

We have shown that integration doesn't work.  Unfortunately, neither does minimization over departure times; that too violates the triangle inequality.  The minima contributing to $T(a,b)$ and $T(b,c)$ could occur at times that are not successive -- and hence there would be no corresponding composable routes.  Suppose that $X$ is the set of $t$ for which $T(a,b,t)=T(a,b)$, the times which determine it via minimization.  Let $Y$ be defined the same way for $T(b,c)$.  There is no guarantee that there exists any $t\in X$ for which $t+T(a,b,t) \in Y$.  

As an example, consider points $a$, $b$, and $c$ with a single route forward between each pair of points. $T(r_1(a,b,t))$ is 10 minutes at 5 PM and 1 hour at all other times, so $T(a,b)= 10m$. Similarly, $T(r_2(b,c,t))$ is 10 minutes at 10 AM and 1 hour at all other times, so $T(b,c)= 10m$.  Finally, $T(r_3(a,c,t))$ is 45 minutes at all times.  Clearly $T(a,c)$ would be 45m, violating the triangle inequality.  At any given time of day, the composed routes $a$ to $b$ followed by $b$ to $c$ yield a combined travel time of 70m.  Therefore, the quickest route from $a$ to $c$ always will be $r_3$ and take 45m.  However, this means that $T(a,c)>T(a,b)+T(b,c)$.  

A min-min solution thus proves unsuitable as well.  However, a max-min approach overcomes these problems\footnote{The existence of a max-min metric should not surprise us.  Similar solutions arise for different reasons in a variety of areas.  Of particular note are the minimax solutions of game theory\cite{minimax}, which minimize the severity of one's maximum loss, and the Hausdorff Distance between non-empty compact subsets of a metric space\cite{hausdorff}. }.  We maximize over both $t\in P$ and the direction.  It is important to note that the two maximization operations commute with one another but \textit{not} with the minimization over routes.

\begin{equation}
T(a,b)= \max_{(a',b')\in \{(a,b),(b,a)\}} \max_{t\in P} \min_{r\in R(a',b',t)} T(r(a',b',t))
\label{eq:finalmetric}
\end{equation}

The problem that disqualified a min-min choice doesn't arise; we no longer require that there exist minimizing routes which are composable.  Returning to our example, $T(a,b)$ would be 1h, $T(b,c)$ would be 1h, and $T(a,c)$ would be 45m.  There exist 70m routes from $a$ to $c$, but at any given time the quickest route from $a$ to c is 45m.  So the maximum over times is 45m.

Clearly $T(a,a)=0$ and $T(a,b)>0$ for $a\ne b$.  Let us first prove the triangle inequality in the unsymmetrized case.  To this end, we define

$$T_U(a,b)= \max_{t\in P} \min_{r\in R(a,b,t)} T(r(a,b,t))$$

and 

$$T(a,b)= \max_{(a',b')\in \{(a,b),(b,a)\}} T_U(a',b')$$

Let us suppose that $T_U$ does not obey the triangle inequality. Then there must exist some time $t'$ for which the shortest travel-time from $a$ to $c$ departing at $t'$ takes longer than the sum of the \textit{worst} shortest travel times from $a$ to $b$ and from $b$ to $c$ with any departure times\footnote{Of course, these may take place with different departure times and may not represent routes which are composable.}.  That is,

$$ T(a,c,t') > \max_t T(a,b,t) + \max_t T(b,c,t)$$

But we also know that 

$$\max_t T(a,b,t)\ge T(a,b,t')$$ 

for any specific $t'\in P$ and 

$$\max_t T(b,c,t)\ge T(a,b,t'+T(a,b,t'))$$

for any $t'+T(a,b,t')\in P$.   For now, let us ignore the boundary question of what happens when the latter condition does not hold and assume it does.  Putting this together we have

$$T(a,c,t')> T(a,b,t') + T(a,b,t'+T(a,b,t'))$$

This directly contradicts equation~\ref{eq:newtriang}, which must hold.  Therefore such a $t'$ cannot exist, and $T_U$ must satisfy the triangle inequality. 

Returning to the symmetrized $T$, it may seem possible that our argument breaks down.  The premise of equation~\ref{eq:newtriang} may be violated if the relevant routes no longer are composable.  They could be in opposite directions, for example.  

However, it is easy to see that the triangle inequality endures.  We need only note that 

$$T_U(a,c)\leq T_U(a,b)+T_U(b,c)$$
$$T_U(c,a)\leq T_U(b,a)+T_U(c,b)$$

It always is true that $\max(x_1+y_1,x_2+y_2)\leq \max(x_1,x_2)+\max(y_1,y_2)$.  So,

$$\max(T_U(a,c),T_U(c,a))\leq \max(T_U(a,b),T_U(b,a))+\max(T_U(b,c),T_U(c,b))$$

and the triangle inequality is preserved.

\section{Regularization of the Boundary}

One advantage of our max-min definition of $T$ is that it is agnostic to the actual set of routes deemed admissible.  Our proof relies on the composability of routes.  As long as a given route is included either in all calculations or in none we are fine.  And as long as our consistency assumptions are adhered to by the routes we \textit{do} include, we are free to exclude or include routes as we desire.  This gives us a flexibility that would be absent from a definition which relies more heavily on the details of the routes involved.  

We have been cavalier about the boundary of the period $P$ over which we maximize.  In particular, our proof of the triangle inequality skirted this issue.   It may seem that the same freedom we just mentioned allows us complete latitude to exclude routes that depart after a certain time.  However, the boundary can lead to subtle violations of equation~\ref{eq:tti} if we are careless.  Note that because all routes are forward-facing in time, only the latter boundary of $P$ is of concern\footnote{If instead we had decided to designate routes by their end-times -- and changed our definitions and equations accordingly -- then the early boundary would be the issue.}.  

An example may help illustrate the difficulties that can arise. For simplicity we use $T_U$ to illustrate, but identical considerations hold for $T$.  Let $P$ be some suitable period.  Suppose we have $3$ points $a$, $b$, and $c$, and that routes and travel times are stationary within $P$.  The travel times are as follows:

\begin{itemize}
\item $T(r_1(a,c,t))= 2h$ for all $t\in P$.
\item $T(r_2(a,b,t))= 30m$ for all $t\in P$.
\item $T(r_3(b,c,t))= 30m$ for all $t\in P$.
\end{itemize}

Up until 30 minutes before the end of $P$, we have $T(a,b,t)=30m$, $T(b,c,t)=30m$ and $T(a,c,t)=1h$ (using the route via $b$).  If we perform the same calculation 15 minutes prior to the end of $P$, the first two are the same but $T(a,c,t)=2h$.  This means $T_U(a,c)=2h$ when we maximize, and the triangle inequality is violated.   

Where did we go wrong?  This should be no different from $r_3$ actually ceasing to exist at a certain time, a legitimate scenario.  Our proof of the triangle inequality is predicated on our consistency assumptions.  In our case, the existence requirement is violated, which in turn invalidates equation~\ref{eq:newtriang}.  Because we restrict ourselves to routes with departure times in $P$, no route exists going from $b$ to $c$ at $t+T(a,b,t)$.  The problem would be the same if we eliminated $r_3$ at any time \textit{within} $P$.  If we remedy this by adding a suitable route $r'(b,c,t+T(a,b,t))$, then either it would make $T_U(b,c)$ larger or, composed with some route $r(a,b,t)$, would provide means for $T(a,c,t)$ to be smaller and thus remove it from contributing to $T_U(a,c)$ during the maximization.  

The problem is the same whether we violate the existence requirement because of the boundary or because of a real paucity of routes.  In the latter case, the fault would lay with the data set, model, or practitioner.  The problem at hand simply would not satisfy the requirements for use of our metric.  But the boundary is present in any realistic situation, so how do we deal with it? 

Extending the boundary accomplishes nothing.  Attempting to iteratively add routes that depart after $P$ or remove ones within $P$ to maintain consistency is difficult and fraught with problems.  Fortunately, there is a simple way to achieve the same effect with little effort.  

The boundary is the problem, so we eliminate it.  Suppose $P=[t_0,t_1]$.  By assumption, there exists at least one route between any pair of points $a$ and $b$ for each departure time $t_0 \leq t \leq t_1$.  We extend this by defining routes $r^*(a,b,t)$ for all locations $a\ne b\in M$ and times $t>t_1$ and such that $T(r^*(a,b,t))=\epsilon$, where $\epsilon$ is very small\footnote{A safe choice would be to make it shorter than the shortest non-trivial route between any points in $M$ and at any time in $P$.  Specifically, we could choose $\epsilon=\frac{1}{2} \min_{a\ne b\in M} \min_{t\in P} \min_{r\in R(a,b,t)} T(r(a,b,t))$.}.  For completeness, we can add the trivial $e(a,a,t)$ routes for $t>t_1$ as well.  The maximization then is over $t\in [t_0,\infty]$ and we have no boundary.  Clearly, $T(a,b,t)=T(r^*(a,b,t))$ for $t>t_1$ and does not contribute to $T_U(a,b)$ at all\footnote{If we later require continuity or differentiability of the travel times $T(a,b,t)$ as functions of $t$, we could instead smoothly continue them past $t_1$ using a bump function.}.  In effect, the regularization selectively eliminates the contribution of certain $T(a,b,t)$ to $T_U(a,b)$ in a way which varies with $a$ and $b$.  

How would this work in our example?  $T_U(a,b)$ and $T_U(b,c)$ are unaffected.  However, near the boundary there now is a route $r_1(a,b,t)+r^*(b,c,t+T(r_1(a,b,t)))$ with travel time $30m+\epsilon$.  This yields an inaccurate $T(a,c,t)$; however, there is no contribution to $T_U(a,c)$ from that value because it is eliminated during the maximization.  What we effectively have done is vary the maximization period $P$ for different pairs of points.  In the case of $(a,c)$, it runs up to 30 minutes before $t_1$.  In the case of $(a,b)$ it will get closer to $t_2$, but how close depends on which routes minimize $T(a,b,t)$ for the various $t$.  Our regularization method has performed the very difficult task of assessing these things for us.  Note that we have done nothing illegal\footnote{The standard refrain of scoundrels, to be sure. Hopefully, we offer a compelling defense.}.  Time-dependent results outside the effective period $P(a,b)$ for each pair will be meaningless -- but the overall $T(a,b)$ is both meaningful and a metric.  

It may be tempting to use the same approach to supply missing routes in real problems which violate the existence requirement.  However, great care should be taken in doing so.  The routes we added all were outside of our domain ($M,P$).  If we add deminimus routes \textit{within} that domain then they may yield tiny $T(a,b)$ values between their endpoints and distort all sorts of other $T$ values -- albeit in a way which preserves the triangle inequality.

\section{Design Decisions}

There are several decisions that must be made when adopting our approach, and they can significantly affect the utility of the results.  

The first question is how we define ``travel'', what modes and routes are admissible.  We may wish to focus on a particular means of travel -- vehicle, public transport, or walking, for example -- or allow arbitrary combinations.   

One could be tempted to modify the definition of travel-time to represent some sort of more general utility.  Perhaps it is desirable to include the cost of travel; reflect quality of life issues such as inconvenience, danger, or discomfort; or incorporate other factors.  Unfortunately, we do not have license to do so.  The travel times are used to determine end-to-end composability of routes.  They have an absolute meaning within the model and cannot be replaced with arbitrary nonnegative quantities.  While conceivably the two uses of $T$ could be separated, doing so is beyond our present scope.  What we should and must include in our travel times are any time-related features:  parking, waiting for buses or trains, etc.  $T(r(a,b,t))$ must represent a true location-to-location travel time, or composition will be meaningless.

Depending on the situation at hand, there may be real routes which start and end within $M$ but pass outside of it.  For example, when traveling within a city it sometimes may be faster to get on a highway, leave the city, and re-enter it near the destination.  We are free to discard any routes that cross the boundaries of the region, as long as we adhere to our consistency conditions.   Whether it is desirable to do so depends on the availability of data for such routes as well as the intended meaning of the metric.  If we exclude them, the travel time between certain points may appear larger than would be experienced.  On the other hand, if we include them then some results may seem inexplicable when viewed entirely in the context of $M$.

It may be important to understand how stable the travel-times are.  If there are many efficient routes in both directions between two points, then it is fair to say that those points are highly accessible to one another.  On the other hand, if there only is one such route then any disruption to it could cause a substantial change in travel time.  Our metric does not distinguish these cases, and doing so may constitute an important supplement.  There are various ways to accomplish this.  As a simple example, one could exclude the fastest route between any two points from each $R(a,b,t)$ before performing our calculations.  As long as the existence requirement remains intact, this is entirely legitimate.  A comparison of the original metric to this modified one could offer a measure of its stability. 
 
Our travel time metric $T(a,b)$ is based on a particular period $P$.  If we regularize the boundary, we do so relative to that period.  $P$ contains the real data from which we derive $T$.  Though we maximize over departure times in $P$, there is no reason that our metric cannot itself be a function of time.  As an example, we could employ a rolling period $P(t)$ defined as a 7 day interval $[t-7d,t]$.  At risk of butchering our notation, we could designate the associated metric $T(P(t),a,b)$.  For each choice of $t$, we have a period $P$ and a metric $T$.  It is important to note the difference between the time-dependence of such a metric and that present in the $T(a,b,t)$ function.  The latter measures the quickest travel between points for a given departure time.  It does not constitute a metric on $M$ if viewed as a function of $a$ and $b$ at fixed $t$.  Rather, the $T(a,b,t)$ at different times are related via equation~\ref{eq:newtriang}.  On the other hand, $T(P(t),a,b)$ \textit{is} a metric on $M$ at each time $t$.  The relationships amongst these $T$ values are contemporaneous.   If we wish to be truly pedantic, we could write our metric with all its major dependencies as $T(M,P,R,a,b)$ where $R$ is a function on $M\times M\times P$ whose value is the set $R(a,b,t)$ for any $a$, $b$, and $t$. 

Finally, our metric ignores any question of capacity -- an important consideration in urban planning.  There are a couple of ways this may be addressed within our framework.  One could compare the metric at different capacity levels.  With this approach, every route is assigned a maximum capacity.  The metric is then computed using various traffic-level scenarios.  In a given scenario, a certain existing volume is assumed along each route.  The metric measures travel-times for a marginal traveler, so any routes which are at capacity are denied the traveler.  As long as we don't end up removing enough routes to violate the existence requirement, we obtain a metric for each scenario.  A comparison of these may yield insight into the dependence of travel-times on traffic-volume.  Alternately, capacity could be accounted for in a modification of travel times along individual routes.  A route is between fixed endpoints, has a fixed starting time, and reflects a particular mode of transport.  It makes perfect sense to adjust the associated travel time to reflect a realistic traffic scenario.  For example, we could assume a certain probability of having to wait for a second bus, or we could extend travel times for drivers.  Care must be taken to do so in a manner consistent with our assumptions, however.

\end{document}